\documentclass[amstex,12 pt]{article} 
\usepackage{amsmath}
\usepackage{amsfonts}
\usepackage{MnSymbol}

\newtheorem{Th}{Theorem}[section]
\newtheorem{Prop}[Th]{Proposition}
\newtheorem{Lem}[Th]{Lemma}

\newenvironment{pf}{\noindent{\bf Proof.}}{\CQFD
}
\newcommand{\CQFD}
{%
\mbox{}%
\nolinebreak%
\hfill%
\rule{2mm}{2mm}%
\medbreak%
\par%
}

\newcommand{\R}{\mathbb{R}}

\newcommand{\C}{\mathbb{C}}

\newcommand{\B}{\mathbb{B}}

\def\Ho{\vbox{\offinterlineskip\hbox{\kern 3pt$\scriptstyle\circ$}
\kern
1pt\hbox{$H$}}}

\begin{document}

\title
{Splitting Parabolic Manifolds}

\author{Morris KALKA and Giorgio PATRIZIO}
\date{ }

\footnotetext{Much of this work was done while Kalka was
visiting the University of Florence and G. Patrizio   Tulane
University. The Authors thank the institutions for their support. G. Patrizio
acknowledges the support of MIUR PRIN 2010-11 ``Variet\`a reali e complesse: geometria, topologia e analisi armonica'' and the collaboration with GNSAGA of INdAM.}

\maketitle

\begin{abstract} 
\noindent We study the geometric properties of complex manifolds possessing a pair of plurisubharmonic functions satisfying Monge-Amp\`ere type of condition.  
The results are applied to characterize complex manifolds biholomorphic to  $\C^{N}$ viewed as a product of lower dimensional complex euclidean spaces.\\
\end{abstract}


\section {Introduction}

A well known theorem of Stoll (\cite{St}, see also \cite{Burns,Patrizio,Wong}) characterizes $\C^{n}$ (in fact the pair ($\C^{n}, \tau_{n}=\Vert \cdot \Vert^{2}$) up to biholomorphic map as the unique unbounded $n$--dimensional strictly parabolic  manifold i.e. $n$--complex manifold $M$ equipped with a $C^{\infty}$ strictly plurisubharmonic exhaustion $\rho : M \to [0,+\infty)$ such that the function $\log \rho$ satisfies the 
complex homogeneous Monge-Amp\`ere equation $(dd^{c} \log \rho)^{n}=0$ on $\rho>0$. Such a $\rho$ is called a {\sl (unbounded) strictly parabolic exhaustion} of $M$ and the pair $(M,\rho)$ {\sl  (unbounded) strictly parabolic  manifold}. Indeed,if $(M,\rho)$ is an $n$--dimensional  (unbounded) strictly parabolic  manifold, there exists a biholomorphic map $F\colon  \C^{n} \to M$ such that 
$\rho(F(Z))=\Vert Z \Vert^{2}=\tau_{n}(Z)$. Notice that here we do not require in the definition of strictly parabolic exhaustion that $\log \rho$ is plurisubharmonic as in  \cite{St}; indeed  Wong in (\cite{Wong}) proved that this property follows from the other assumptions.

In this paper we are interested in the characterization of $(\C^{N},\tau_{N})$ as product of unbounded strictly parabolic manifolds i.e. of $(\C^{k},\tau_{k})$ and $(\C^{l},\tau_{l})$ for $k,l>0$ with $N=k+l$. More precisely we like to look at $\C^{N}$ and characterize it from the following point of view.
\\
For $N\geq 2$ and $k,l>0$ with $k+l=N$, let $\eta_{k}\colon \C^{N} \to [0,+\infty)$ and  $\eta_{l}\colon \C^{N} \to [0,+\infty)$ be defined for $Z=(z_{1},\dots,z_{k},z_{k+1},\dots z_{N})$ respectively by

\begin{align} \eta_{k}(Z)= \vert z_{1}\vert^{2}+ \dots\vert z_{k}\vert^{2}
\hskip0,5cm
{\rm and} \hskip0,5cm
 \eta_{l}(Z)= \vert z_{k+1}\vert^{2}+ \dots\vert z_{N}\vert^{2}. 
\end{align}
\\
One checks easily that  the following properties hold true:
\begin{itemize}
\item[i)] $\eta_{k},\eta_{l}\in {\cal PSH}(\C^{N})\cap C^{\infty}(\C^{N})$,
\item[ii)] $(dd^{c}\log \eta_{k})^{k}=0$ on $\{\eta_{k}>0\}$ \hskip0,3cm and \hskip0,3cm $(dd^{c}\log \eta_{l})^{l}=0$ on $\{\eta_{l}>0\} $,
\item[iii)] $(dd^{c}\eta_{k})^{k}\wedge (dd^{c}\eta_{l})^{l}>0$ on $\C^{N}$,
\item[iv)] $\eta_{k}+\eta_{l}$ is an exhaustion of $\C^{N}$.
\end{itemize}
It is also straightforward to verify that 
\begin{itemize}
\item[v)] $(dd^{c}\eta_{k})^{k}\wedge d\eta_{k}\wedge d^{c}\eta_{k}=0$ \hskip0,3cm and \hskip0,3cm 
 $(dd^{c}\eta_{k})^{k-1}\wedge d\eta_{k}\wedge d^{c}\eta_{k}\neq0$,
\item[v')] $(dd^{c}\eta_{l})^{l}\wedge d\eta_{l}\wedge d^{c}\eta_{l}=0$ \hskip0,3cm and \hskip0,3cm 
 $(dd^{c}\eta_{l})^{l-1}\wedge d\eta_{l}\wedge d^{c}\eta_{l}\neq0$
\end{itemize}
and that each level set of $\eta_{k}$ and $\eta_{l}$ are foliated respectively by closed $l$-dimensional submanifolds (in fact $l$-complex affine subspaces) and by closed $k$-dimensional submanifolds (in fact $k$-complex affine subspaces). Indeed the level sets of $\eta_{k}$ and $\eta_{l}$ are respectively foliated by the $l$-complex affine subspaces tangent to the annihilator of $dd^{c}\eta_{k}$ and the $k$-complex affine subspaces tangent to the annihilator of $dd^{c}\eta_{l}$. 
While all other level sets  of $\eta_{k}$ and $\eta_{l}$ are Levi flat smooth real hypersurfaces foliated by the leaves of the $u$ foliation and by the leaves of the $v$ foliation, the minimal sets $\{\eta_{k}=0\}$ and $\{\eta_{l}=0\}$ are exceptional: they consists only respectively of one $u$ leaf and of one $v$ leaf
and therefore, in particular are closed complex analytic submanifolds of $\C^{N}$.
\\ \\
We like to show how these properties may be used to characterize up to biholomorphic maps the triple
$(\C^{N},\eta_{k},\eta_{l})$. In order to state our main result, we need the following:
\\ \\
{\bf Definition.} Let $M$ be a complex manifold of dimension $N=k+l\geq2$ for some $k,l>0$ and 
$u,v : M \to [0,+\infty)$. We say that the triple  $(M,u,v)$ is a $(k,l)${\sl --splitting parabolic manifold} if the following conditions hold:
 \begin{itemize}
\item[(A1)] $u,v\in {\cal PSH}(M)\cap C^{\infty}(M)$,
\item[(A2)] $(dd^{c}\log u)^{k}=0$ on $\{u>0\}$ \hskip0,3cm and \hskip0,3cm $(dd^{c}\log v)^{l}=0$ on $\{v>0\} $,
\item[(A3)] $(dd^{c}u)^{k}\wedge (dd^{c}v)^{l}>0$ on $M$,
\item[(A4)] $u+v$ is an exhaustion of $M$.
\end{itemize}
Using the theory of complex Monge-Amp\`ere foliations, it is possible to reconstruct the geometry of splitting parabolic manifolds. In particular we  show that there exist foliations ${\cal F}_{u}$ and  ${\cal F}_{v}$ of $M$ respectively by  $l$-dimensional submanifolds  and distribution ${\rm Ann} (dd^{c}u)^{k}$ and by closed $k$-dimensional submanifolds and distribution ${\rm Ann} (dd^{c}v)^{l}$ such that  
\begin{itemize}
\item[--] {\sl $u$ is constant along each leaf of ${\cal F}_{u}$ and  ${\cal F}_{u}$ restricts to a foliation of any level set of $u$;}
\item[--] {\sl $v$  is constant along each leaf of ${\cal F}_{v}$ and  ${\cal F}_{v}$ restricts to a foliation of any level set of $v$.}
\end{itemize} 
Under additional global assumptions on these foliations, more can be proved and it turns out that for
$(k,l)$--splitting parabolic manifolds one may reconstruct almost entirely  
the geometric picture of $\C^{N}$ as product lower dimensional complex euclidean spaces. In fact, assuming that the ${\cal F}_{u}$ and  ${\cal F}_{v}$ foliations have closed leaves, one shows that:

\begin{itemize}
\item[--] {\sl each leaf of the ${\cal L}_{u}$ foliation is biholomorphic to $\C^{l}$ and  each leaf of the ${\cal L}_{v}$
foliation is biholomorphic to $\C^{k}$ and both the ${\cal L}_{u}$ foliation and the ${\cal L}_{v}$
foliation are holomorphic.}
\item[--] {\sl the set $L^{u}_{0}=\{u=0\}$ is complex analytic of dimension $l$ and it is a leaf of the ${\cal L}_{u}$ foliation 
as well as the set $L^{v}_{0}=\{v=0\}$ is complex analytic of dimension $k$ and it is a leaf of the ${\cal L}_{v}$ foliation;}
\item[--] {\sl for every leaf $L$ of the ${\cal L}_{u}$ foliation, there exists a unique point $O_{L}$ such that 
$L\cap L^{v}_{0}=\{O_{L}\}$ and for every leaf $L'$ of the ${\cal L}_{v}$ foliation, there exists a unique point $O_{L'}$ such that 
$L'\cap L^{u}_{0}=\{O_{L'}\}$}
\end{itemize}

The existence and the properties of the foliations ${\cal F}_{u}$, ${\cal F}_{v}$  under the global topology assumption on their topology, namely the closeness of the leaves, allows to characterize $\C^{N}$ with this kind of product structure:

\vskip0,3cm
\noindent
{\bf Characterization Theorem.} 
{\sl Let $M$ be a complex manifold of dimension $N=k+l\geq 2$ for some $k,l>0$ and 
$u,v : M \to [0,+\infty)$ such that the triple $(M,u,v)$ is a $(k,l)$--splitting parabolic manifold. Then
the leaves of the foliations ${\cal F}_{u}$ and  ${\cal F}_{v}$ are closed if and only if there exists a biholomorphic map 
$F : \C^{N} \to M$ such that for $Z\in \C^{N}$ one has
$u\circ F(Z)= \eta_{k}(Z)$  and 
 $u\circ F(Z)= \eta_{l}(Z).$ }

\section {Geometry of Splitting Parabolic Manifolds.}

We now describe the geometry of splitting parabolic manifolds and provide the steps necessary to prove our main result. We start with two technical remarks which sprout directly from the assumptions.

\begin{Lem} \label{compute} Let $M$ be a complex manifold of dimension $N$ with $N=k+l$ for some $k,l>0$.
Suppose that there exist functions $u, v: M \to [0,+\infty)$ such that  $(A1), (A2), (A3)$ hold. Then:
\begin{align} \label{leviu} 
\left\{
\begin{array}{l}
(dd^{c}u)^{k}\wedge du\wedge d^{c}u=0 \; on\; M, \\ \\ 
 (dd^{c}u)^{k-1}\wedge du\wedge d^{c}u\neq0 \;on\; \{u>0\}, 
 \end{array}
 \right.
 \end{align}
\begin{align} \label{leviv}
 \left\{
\begin{array}{l}
(dd^{c}v)^{l}\wedge dv\wedge d^{c}v=0 \; on\; M,
\\ \\
 (dd^{c}v)^{l-1}\wedge dv\wedge d^{c}v\neq0 \;on\; \{v>0\};
  \end{array}
 \right.
 \end{align}
 \begin{align} \label{mau}(dd^{c}u)^{k+1}=0 \hskip0,3cm with \hskip0,3cm 
 (dd^{c}u)^{k}\neq0 \; on\; M;
 \end{align}
\begin{align} \label{mav}(dd^{c}v)^{l+1}=0  \hskip0,3cm with \hskip0,3cm 
 (dd^{c}v)^{l}\neq0 \; on\; M;
 \end{align}
 \begin{align} \label{level} du\neq0  \;on\; \{u>0\}, \hskip0,3cm with \hskip0,3cm 
 dv\neq0 \;on\; \{v>0\}.
 \end{align}

\end{Lem}

\begin{pf} One has immediately:
\begin{align} \label{derlog}u dd^{c} u = u^{2} dd^{c}\log u + du\wedge d^{c}u,
 \end{align}
 so that, taking exterior powers and using $(A2)$, one has on $\{u>0\}$:
 \begin{align}\label{kpower} u^{k}(dd^{c}u)^{k}= k u^{2(k-1)}(dd^{c}\log u)^{k-1}\wedge du\wedge d^{c}u. \end{align}
 On $\{u=0\}$, the minimal set of $u$, we have $du=0$ so that
 $(dd^{c}u)^{k}\wedge du\wedge d^{c}u=0$. 
 On the other hand $(dd^{c}u)^{k}\wedge du\wedge d^{c}u=0$ on $\{u>0\}$ is an immediate consequence of (\ref{kpower}).  Furtermore from $(A3)$ it follows that $(dd^{c}u)^{k}\neq0$, so that , from (\ref{derlog}) and (\ref{kpower}), we get 
  \begin{align}\label{k1power} 
  \begin{array}{ll} u^{k-1}(dd^{c}u)^{k-1}\wedge du\wedge d^{c}u &= 
  (u^{2}dd^{c}\log u + du\wedge d^{c} u)^{k-1}\wedge du\wedge d^{c}u \\ & \\ &
  = u^{2(k-1)}(dd^{c}\log u)^{k-1}\wedge du\wedge d^{c}u.  \\ & \\ &
  =\frac{1}{k} u^{k} (dd^{c}u)^{k}\neq0
  \end{array}
  \end{align}
so that we have proved (\ref{leviu}). The same applied for $v$, gives (\ref{leviv}).
From (\ref{derlog}), taking the $(k+1)$--exterior power one has 
$(dd^{c}u)^{k+1}=0$ while, as we observed, from $(A3)$ it follows
 $(dd^{c}u)^{k}\neq0$ so that (\ref{mau}) holds. The same argument
 shows (\ref{mav}). Finally (\ref{level}) is immediate from the second part of (\ref{leviu}) and (\ref{leviv}).

\end{pf}

\begin{Lem} Let $M$ be a complex manifold of dimension $N$ with $N=k+l$ for some $k,l>0$.
Suppose that there exist functions $u, v: M \to [0,+\infty)$ such that  $(A1), (A2), (A3), (A4)$ hold. Then $\tau=u+v$ is a strictly plurisubharmonic exhaustion of $M$ which, therefore, is a Stein manifold.
\end{Lem}

\begin{pf}  Since $u$ and $v$ are plurisubharmonic, we have $dd^{c} \tau \geq 0$ on $M$. It is therefore enough to show that $(dd^{c} \tau)^{N} \neq 0$. This, using  (\ref{mau}),   (\ref{mav})
and $ (A3)$, is consequence of the following:
\begin{align}
\begin{array}{ll} (dd^{c}\tau)^{N} &= (dd^{c}u+dd^{c}v)^{N}  
=\sum_{j=0}^{N} {N \choose j} (dd^{c}u)^{j}\wedge (dd^{c}v)^{N-j} \\ & \\ &
={N \choose k} (dd^{c}u)^{k}\wedge (dd^{c}v)^{l}>0.
  \end{array}
  \end{align}
\end{pf}

The main geometric feature of a $(k,l)$--splitting parabolic manifold $(M,u,v)$ is the existence of foliations associated to the functions $u$ and $v$ which can be proved even under milder conditions:

\begin{Prop} Let $M$ be a complex manifold of dimension $N$ with $N=k+l$ for some $k,l>0$.
Suppose that there exist functions $u, v: M \to [0,+\infty)$ such that  $(A1), (A2), (A3)$ hold. Then there exist two foliations ${\cal L}_{u}$, ${\cal L}_{v}$ of $M$ with the following properties:
\begin{itemize}
\item[(i)] ${\cal L}_{u}$ is a foliation in complex $l$-dimensional complex manifolds whose associated distribution is given by the annihilator of the form $dd^{c}u$. Furthemore the function $u$ is constant along each leaf of ${\cal L}_{u}$ and, for all $t>0$, ${\cal L}_{u}$ defines a foliation in $l$-dimensional complex manifolds of the Levi flat real hypersurface 
$\{u=t\}$.
\item[(ii)] ${\cal L}_{v}$ is a foliation in complex $k$-dimensional complex manifolds whose associated distribution is given by the annihilator of the form $dd^{c}v$. Furthemore the function $v$ is constant along each leaf of ${\cal L}_{v}$ and, for all $s>0$, ${\cal L}_{v}$ defines a foliation in $k$-dimensional complex manifolds of the Levi flat real hypersurface 
$\{v=s\}$.
\end{itemize}
\end{Prop}

\begin{pf}  Under the assumptions  $(A1), (A2), (A3)$ we have the conclusions of Lemma
\ref{compute}. It is well known (see Theorem 2.4 of \cite{Bedford-Kalka} for instance) that, as a consequence of (\ref{mau}), the annihilator of the form $dd^{c}u$ is of (complex) rank $l$ and defines a 
foliation ${\cal L}_{u}$ in $l$-dimensional complex submanifolds and, as a consequence of (\ref{mav}),
 the annihilator of the form $dd^{c}v$ is of (complex) rank $k$ and defines a 
foliation ${\cal L}_{v}$ in $k$-dimensional complex submanifolds. On the other hand, since  (\ref{leviu}) holds on $ \{u>0\}$ and (\ref{leviu}) holds on $ \{v>0\}$, from  Theorem 5.3 and Corollary 5.4 of \cite{Bedford-Kalka}
we get that the function $u$ is constant along each leaf of ${\cal L}_{u}$ and, for all $t>0$, ${\cal L}_{u}$ defines a foliation in $l$-dimensional complex manifolds of the Levi flat real hypersurface 
$\{u=t\}$ and that the function $v$ is constant along each leaf of ${\cal L}_{v}$ and, for all $s>0$, ${\cal L}_{v}$ defines a foliation in $k$-dimensional complex manifolds of the Levi flat real hypersurface 
$\{v=s\}$.
\end{pf}

Under the more stringent topological requirement that the leaves of the ${\cal L}_{u}$, ${\cal L}_{v}$ are closed, it can be shown that the leaves are all parabolic and that, consequently, the foliations are holomorphic:

\begin{Th}  \label{parabolicleaves} Let $M$ be a complex manifold of dimension $N=k+l$ for some $k,l>0$ and 
$u,v : M \to [0,+\infty)$ such that the triple $(M,u,v)$ is a $(k,l)$--splitting parabolic manifold such that the leaves of the foliations ${\cal L}_{u}$, ${\cal L}_{v}$ of $M$ are closed. Then:
\begin{itemize}
\item[(i)] If  $L$ is a leaf of the ${\cal L}_{u}$ foliation, then  $(L,v)$ is an $l$--dimensional unbounded strictly parabolic manifold and if $L'$ is a leaf of the ${\cal L}_{v}$ foliation,  then  $(L',u)$ is a $k$--dimensional unbounded strictly parabolic manifold; 
\item[(ii)]  each leaf of the ${\cal L}_{u}$ foliation is biholomorphic to $\C^{l}$ and  each leaf of the ${\cal L}_{v}$
foliation is biholomorphic to $\C^{k}$;
\item[(iii)]  the ${\cal L}_{u}$ foliation and the ${\cal L}_{v}$
foliation are holomorphic.
\end{itemize}
\end{Th}
\begin{pf} $(i)$: If  $L$ is a leaf of the ${\cal L}_{u}$ foliation, since $L$ is closed, the restriction of the exhaustion $u+v$ to $L$ is an exhaustion of $L$. On the other hand $L$ is contained in a level set of $u$, hence the restriction of $v$ to $L$ is an exhaustion of $L$. On the other hand the directions tangent to a leaf  of the ${\cal L}_{u}$ foliation are in the annihilator of $dd^{c}u$ so that in order to have assumption (A3) satisfied, the plurisubharmonic function $v$ must have strictly plurisubharmonic restriction on each leaf of the ${\cal L}_{u}$ foliation. Finally since (A2) holds, it follows that the restriction of $v$ is a strictly parabolic exhaustion of every leaf of  the ${\cal L}_{u}$ foliation. 
The same argument works for the case of the leaves ${\cal L}_{v}$ foliation.
\\
$(ii)$: Because of $(i)$ Stoll's theorem holds for on for each leaves of both ${\cal L}_{u}$ foliation and  each leaf of the ${\cal L}_{v}$ foliation.
\\
$(iii)$: According to $(ii)$ all the leaves of both ${\cal L}_{u}$ foliation and  each leaf of the ${\cal L}_{v}$ foliationare biholomorphic to a complex euclidean space. The conclusion then is direct consequence of  \cite{Kalka-Patrizio4}  where it is shown that in any codimension parabolic locally Monge-Amp\`ere foliation are holomorphic.  
\end{pf}

The final piece of information needed to proceed  is the observation that the minimal sets of $u$ and $v$ are complex analytic and are leaves of the  ${\cal L}_{u}$ foliation and of the ${\cal L}_{}$ foliation respectively:

\begin{Prop} \label{minimal} Let $M$ be a complex manifold of dimension $N=k+l$ for some $k,l>0$ and 
$u,v : M \to [0,+\infty)$ such that the triple $(M,u,v)$ is a $(k,l)$--splitting parabolic manifold such that the leaves of the foliations ${\cal L}_{u}$, ${\cal L}_{v}$ of $M$ are closed. Then:
\begin{itemize}
\item[(i)] the set $L^{u}_{0}=\{u=0\}$ is complex analytic of dimension $l$ and it is a leaf of the ${\cal L}_{u}$ foliation; 
\item[(ii)]   the set $L^{v}_{0}=\{v=0\}$ is complex analytic of dimension $k$ and it is a leaf of the ${\cal L}_{v}$ foliation.
\item[(iii)] For every leaf $L$ of the ${\cal L}_{u}$ foliation, there exists a unique point $O_{L}$ such that 
$L\cap L^{v}_{0}=\{O_{L}\}$; 
\item[(iv)]   For every leaf $L'$ of the ${\cal L}_{v}$ foliation, there exists a unique point $O_{L'}$ such that 
$L'\cap L^{u}_{0}=\{O_{L'}\}$;
\end{itemize}
\end{Prop}

\begin{pf}  We start proving $(i)$.  The set $L^{u}_{0}=\{u=0\}$  is closed and, since the foliation is holomorphic, for each $p\in L^{u}_{0}$ there exists a neighborhood $U$ and  coordinates $z_{1},\dots,z_{k},z_{k+1},\dots,z_{k+l}$ on $U$ so that the intersection of a leaf with $U$ is given by $z_{1}=c_{1},\dots,z_{k}=c_{k}$
for suitable constants $c_{1},\dots,c_{l}$  and $z_{k+1},\dots,z_{k+l}$ are leaf coordinates. In these coordinates one has that on $U$ the Levi matrix of $u$ is as follows:
\begin{align}\label{Leviatzeroset}\left({u_{r\bar s}}\right)_{U} =
\left(\begin{array}{cccccc}
 0 & \dots& 0 & 0& \dots& 0\\
  \vdots &  \vdots& \vdots&\dots & \dots & \vdots \\ 
  \vdots &  \vdots& 0&0 & \dots &0 \\ 
  0  & \dots  & 0 & u_{{k+1}\, \overline{k+1}}& \dots& u_{{k+1}\, \bar k+l}\\
\vdots &  \vdots& \vdots&\vdots & \vdots & \vdots \\
  0  & \dots  & 0 & u_{k+l\, \overline{k+1}}& \dots& u_{k+l \bar k+l}    
\end{array}\right).
\end{align}
from which it follows that the functions $u_{1},\dots,u_{k}$ are holomorphic on $U$. Since $u_{1}=\dots=u_{k}=0$ along $L^{u}_{0}\cap U$, it follows that $L^{u}_{0}$ is an analytic set of dimension $l$. Since every point of 
$L^{u}_{0}$ is contained in a leaf of the ${\cal L}_{u}$, the analytic set $L^{u}_{0}$ is indeed a leaf.
 The same argument works for   $L^{v}_{0}=\{v=0\}$ and therefore $(ii)$ holds. 
 \par
As for $(iii)$ and $(iv)$, they are consequence of the fact that $u$ and $v$ are strictly parabolic exhaustions respectively
for each leaf  $L'$ of the ${\cal L}_{v}$ foliation and of each leaf $L$ of the ${\cal L}_{u}$ foliation. It is known that a key ingredient of the proof of Stoll's Theorem is the fact that the zero set of a strictly parabolic exhaustion reduces to one single point (see Theorem 2.5 of \cite{St}). This implies both  $(iii)$ and $(iv)$.
\end{pf}

\section {A characterization theorem and final remarks.}

For $N\geq 2$ and $k,l>0$ with $k+l=N$, recall that $\eta_{k}\colon \C^{N} \to \R$ and  $\eta_{l}\colon \C^{N} \to \R$ are the functions defined 
 for $Z=(z_{1},\dots,z_{k},z_{k+1},\dots z_{N})$ respectively by

\begin{align} \eta_{k}(Z)= \vert z_{1}\vert^{2}+ \dots\vert z_{k}\vert^{2}
\hskip0,5cm
{\rm and} \hskip0,5cm
 \eta_{l}(Z)= \vert z_{k+1}\vert^{2}+ \dots\vert z_{N}\vert^{2}. \label{uno}
\end{align}

With these notations, we  now state our characterization result:  

\begin{Th} \label{mainth}
 Let $M$ be a complex manifold of dimension $N=k+l\geq 2$ for some $k,l>0$ and 
$u,v : M \to [0,+\infty)$ such that the triple $(M,u,v)$ is a $(k,l)$--splitting parabolic manifold. Then
the leaves of the foliations ${\cal F}_{u}$ and  ${\cal F}_{v}$ are closed if and only if there exists a biholomorphic map 
$F : \C^{N} \to M$ such that for $Z=(z_{1},\dots,z_{k},z_{k+1},\dots z_{N})\in \C^{N}$ one has
\begin{align} u\circ F(Z)= \eta_{k}(Z) \hskip0,5cm
{\rm and} \hskip0,5cm
 u\circ F(Z)= \eta_{l}(Z). \label{tre}
\end{align}
\end{Th}

\begin{pf} In one direction the result is obvious: if for a $(k,l)$--splitting parabolic manifold $(M,u,v)$
there exists $F : \C^{N} \to M$ such that (\ref{tre})  holds, then $F$ maps the leaves of the foliation of $\C^{N}$ associated to $\eta_{k}$ and the leaves of the foliation of $\C^{N}$ associated to $\eta_{l}$ respectively onto the leaves of the ${\cal L}_{u}$ foliation and the the ${\cal L}_{v}$ foliation of $M$ and therefore the leaves of these are closed. 
\vskip0,05cm
Suppose now that $(M,u,v)$ is a $(k,l)$--splitting parabolic manifold such that
the leaves of the foliations ${\cal F}_{u}$ and  ${\cal F}_{v}$ are all closed. Then the conclusions of Proposition
 \ref{parabolicleaves} and of Proposition \ref{minimal} hold true and we shall use the notations introduced there. Futhermore we set $L^{u}_{0}\cap L^{v}_{0}=\{\bf O\}$ and, for any $p\in M$, denote  by 
 $L^{u}(p)$ and  $L^{v}(p)$  the unique leaves respectively of the ${\cal F}_{u}$ foliation and  of the ${\cal F}_{v}$ foliation passing through $p$. Notice that with that, in this case we have $L^{u}({\bf O})=L^{u}_{0}$ and  $L^{v}({\bf O})=L^{v}_{0}$
 and we shall insist with the notations $L^{u}_{0}$ and $L^{v}_{0}$ for these special leaves.  
 Thus, there exist maps of class $C^{\infty}$
\begin{align} \label{globalparam}
\Phi^{u}\colon \C^{l}\times L^{v}_{0}\to M \hskip0,2cm{\rm and} \hskip0,2cm
 \Phi^{v}\colon \C^{k}\times L^{u}_{0}\to M\end{align}
 such that for all $x\in  L^{u}_{0}$ and for all $y\in  L^{v}_{0}$ the maps
\begin{align} \label{leafparam}
 \Phi^{v}(\cdot,x)\colon \C^{k}\to  L^{v}(x) \hskip0,1cm{\rm and} \hskip0,1cm
 \Phi^{u}(\cdot,y)\colon \C^{l}\to  L^{u}(y)
\end{align}
are biholomorphic and such that for all $x\in  L^{u}_{0}$ and for all $y\in  L^{v}_{0}$
\begin{align} \label{leafisom}
v\circ\Phi^{u}(w_{1}\dots,w_{l},y)= \vert w_{1}\vert^{2}+ \dots\vert w_{l}\vert^{2}
\hskip0,2cm{\rm and} \hskip0,2cm
u\circ\Phi^{v}(z_{1}\dots,z_{k},x)= \vert z_{1}\vert^{2}+ \dots\vert z_{k}\vert^{2}.
\end{align}
The existence of the maps $\Phi^{u}$ and $\Phi^{v}$ satisfying (\ref{leafparam}) and (\ref{leafisom}) follows exactly from the content of Stoll's characterization 
(\cite{St}) of unbounded strictly parabolic manifolds. The smoothness of the maps defined in 
(\ref{globalparam}) is consequence of the construction needed to prove Stoll's theorem. In fact, once the 
the complex euclidean space is identified with  the tangent space at the center of the manifold,  the map that
provides the biholomorphism  is exactly the exponential map of the K\"ahler metric whose potential is the strictly parabolic exhaustion (see \cite{St}, \cite{Burns} and \cite{Patrizio}). In our case the stricly parabolic exhaustions are the restrictions of $v$ to the leaves of the ${\cal F}_{u}$ foliation and of $u$ to the leaves of the  ${\cal F}_{v}$ foliation. As $u$ and $v$ are of class $C^{\infty}$ then the smoothness of $\Phi^{u}$ and $\Phi^{v}$ follows as a consequence. 
We need some piece of notation. 
Define projection maps
\begin{align}
\Pi^{u}\colon M\to L^{u}_{0} \hskip1cm {\rm and} \hskip1cm \Pi^{v}\colon M\to L^{v}_{0}
\end{align}
by imposing for $p\in M$:
\begin{align}
\Pi^{u}(p)= L^{v}(p)\cap L^{u}_{0} \hskip1cm {\rm and} \hskip1cm \Pi^{v}(p)= L^{u}(p)\cap L^{v}_{0}.
\end{align}
Since the 
${\cal F}_{u}$ and the  ${\cal F}_{v}$ foliations are holomorphic, the maps $\Pi^{u}$ and $\Pi^{v}$, which are well defined by Proposition \ref{minimal}, are also holomorphic. Finally we denote by $\Phi^{u}_{0}=\Phi^{u}(\cdot,{\bf O})$ and 
$\Phi^{v}_{0}=\Phi^{v}(\cdot,{\bf O})$ the parametrization of the special leaves $L^{u}_{0}$ and $L^{v}_{0}$.
We can now proceed and define a map $G\colon M \to \C^{k+l}=\C^{N}$ that turns out to be the inverse of the required map $F\colon \C^{N} \to M$. We set
\begin{align}
G(p)= \left((\Phi^{v}_{0})^{-1}\circ \Pi^{v}(p) , (\Phi^{u}_{0})^{-1}\circ \Pi^{u}(p) \right) \in \C^{k+l}=\C^{N}
\end{align}
By construction, $G$ is bijective and holomorphic since its $k$--component and  $l$--component are composition of 
holomorphic maps. It follows that $G$ is biholomorphic. Furthermore, again by construction and using the fact that the map of Stoll's theorem preserves strictly parabolic exhaustions, if 
$G(p)=(z,w)$, then $u(p)= \Vert z\Vert_{k}^{2}=\eta_{k}(z,w)$ and $v(p)= \Vert w\Vert_{l}^{2}=\eta_{l}(z,w)$ where
$\Vert \cdot \Vert_{k}$ and $\Vert \cdot \Vert_{l}$ denote respectively the norms of $\C^{k}$ and $\C^{l}$. But then $F=G^{-1}$
has all the properties required by the claim of the Theorem.
\end{pf}
\vskip0,2cm
\noindent
{\bf Remark 1}
Using the notation $Z=(z,w)=(z_{1}\dots,z_{k},w_{1}\dots,w_{l})\in \C^{k+l}=\C^{N}$, a--posteriori, we have that the map $F\colon \C^{N} \to M$ satisfies the following:
\begin{align} \label{mapF}
F(Z)=F(z,w)= \Phi^{v}(z,\Phi^{u}(w,{\bf O}))= \Phi^{u}(w,\Phi^{v}(z,{\bf O}))
\end{align}
where the second equality  it is equivalent to the non trivial fact that the flows associated to the ${\cal F}_{u}$ and to  the  ${\cal F}_{v}$ foliations do commute. 
\vskip0,2cm
\noindent
{\bf Remark 2} As a consequence of Theorem \ref{mainth}, it follows that if  $(M,u,v)$ is a $(k,l)$--splitting parabolic manifold such that the leaves of the foliations ${\cal F}_{u}$ and  ${\cal F}_{v}$ are closed then $u+v$ is a strictly parabolic
exhaustion for $M$ and in particular that $\left(dd^{c}(u+v)\right)^{N}=0$ on $M\setminus\{{\bf O}\}$.  Because of the high non linearity of the equation involved, this cannot be derived directly from the assumptions (A2) and (A3).
\vskip0,3cm
 A key element of the proof of Theorem \ref{mainth} is the fact that if the leaves of the foliations ${\cal F}_{u}$ and  ${\cal F}_{v}$ are closed then the leaves are unbounded strictly parabolic manifolds and hence biholomorphic to complex euclidean spaces. This, in turn, because of the results of 
\cite{Kalka-Patrizio4}, implies that the  the foliations ${\cal F}_{u}$ and  ${\cal F}_{v}$ are holomorphic.
There is no corresponding result for Monge-Amp\`ere foliation with hyperbolic leaves: in fact most 
Monge-Amp\`ere foliation with hyperbolic leaves are not holomorphic. Because of this there is no hope 
for straightforward extensions of the results of this paper to the characterization of products of balls 
\begin{align}
\B^{k}(R_{1})\times\B^{l}(R_{2})= 
 \left\{
(z,w) \in \C^{k+l}=\C^{N} 
\mid  \Vert z\Vert_{k}^{2}<R_{1}^{2},   \Vert w\Vert_{l}^{2}<R_{2}^{2}
\right\}
\end{align}
of complex euclidean spaces. The key to get such extension, in the proof of  the holomorphicity of the foliations one should replace the arguments based on parabolicity of leaves with properties of the foliation along the minimal set of the functions $u$ and $v$. We hope to come back to this problem later on. 
With the techniques developed in this paper and some some obvious adaptations, it can be shown only the following much weaker result:

\begin{Prop}
Let $M$ be a complex manifold of dimension $N=k+l\geq 2$ for some $k,l>0$ and 
$u : M \to [0,R_{1}^{2})$ and  $v : M \to [0,R_{2}^{2})$ functions such that  $(A1), (A2), (A3), (A4)$ hold. If  the leaves of the foliations ${\cal F}_{u}$ and  ${\cal F}_{v}$ are closed and the foliations ${\cal F}_{u}$ and  ${\cal F}_{v}$ are holomorphic, then there exists a biholomorphic map 
$F : \B^{k}(R_{1})\times\B^{l}(R_{2}) \to M$ such that for $Z\in \B^{k}(R_{1})\times\B^{l}(R_{2})$ one has
$u\circ F(Z)= \eta_{k}(Z)$  and 
 $u\circ F(Z)= \eta_{l}(Z).$ 
 \end{Prop}

\bigskip
\bigskip
\font\smallsmc = cmcsc8
\font\smalltt = cmtt8
\font\smallit = cmti8
\hbox{\parindent=0pt\parskip=0pt
\vbox{\baselineskip 9.5 pt \hsize=3.1truein
\obeylines
{\smallsmc
Morris Kalka
Mathematics Department
Tulane University
6823 St. Charles Ave.
New Orleans, LA 70118
USA
}\medskip
{\smallit E-mail}\/: {\smalltt kalka@math.tulane.edu
}
}
\hskip 0.0truecm
\vbox{\baselineskip 9.5 pt \hsize=3.7truein
\obeylines
{\smallsmc
Giorgio Patrizio
Dip. Matematica e Informatica``U. Dini''
Universit\`a di Firenze
Viale Morgani 67/a
I-50134 Firenze
ITALY
}\medskip
{\smallit E-mail}\/: {\smalltt patrizio@unifi.it}
}
}

\end{document}